\documentclass[11pt,a4paper,psamsfonts]{amsart}
\usepackage{amssymb,amscd,amsxtra,calc}
\usepackage{cmmib57}

\theoremstyle{plain}
    \newtheorem{thm}{Theorem}[section] %
    \renewcommand{\thethm}%
    {\arabic{section}.\arabic{thm}}

    \newtheorem{lemma}[thm]{Lemma}

    \newtheorem{theorem}[thm]{Theorem}

\theoremstyle{definition}

    \newtheorem{remark}[thm]{Remark}
\theoremstyle{remark}

    \newtheorem{setup}[thm]{}

\newcommand{\BCC}{\mathbb{C}}

\newcommand{\BPP}{\mathbb{P}}
\newcommand{\BQQ}{\mathbb{Q}}
\newcommand{\BRR}{\mathbb{R}}
\newcommand{\BZZ}{\mathbb{Z}}

\newcommand{\SO}{\mathcal{O}}

\newcommand{\alb}{\operatorname{alb}}
\newcommand{\Aut}{\operatorname{Aut}}

\newcommand{\Bs}{\operatorname{Bs}}

\newcommand{\group}{\operatorname{group}}

\newcommand{\id}{\operatorname{id}}

\newcommand{\NS}{\operatorname{NS}}
\newcommand{\ord}{\operatorname{ord}}
\newcommand{\Per}{\operatorname{Per}}

\newcommand{\Pic}{\operatorname{Pic}}

\newcommand{\Sing}{\operatorname{Sing}}
\newcommand{\SL}{\operatorname{SL}}
\newcommand{\Tr}{\operatorname{Tr}}
\newcommand{\variety}{\operatorname{variety}}

\newcommand{\Alb}{\operatorname{Alb}} 

\newcommand{\ratmap}
{{\,\cdot\negmedspace\cdot\negmedspace\cdot\negmedspace\to\,}}

\begin{document}

\title[The periodic subvarieties]{%
The $g$-periodic subvarieties for an automorphism $g$ of positive entropy on a compact K\"ahler manifold}
\author{De-Qi Zhang}
\address
{%
\textsc{Department of Mathematics} \endgraf
\textsc{National University of Singapore, 2 Science Drive 2,
Singapore 117543}}
\email{matzdq@nus.edu.sg}

\begin{abstract}
For a compact K\"ahler manifold $X$ and a strongly primitive automorphism $g$
of positive entropy, it is shown that
$X$ has at most $\rho(X)$ of
$g$-periodic prime divisors.
When $X$ is a projective threefold, every prime divisor
containing infinitely many $g$-periodic curves, is shown to be $g$-periodic
(a result in the spirit of the Dynamic Manin-Mumford conjecture as in \cite{Zs}).

\end{abstract}


\thanks{The author is supported by an ARF of NUS}

\maketitle


\section{Introduction}

We work over the field $\BCC$ of complex numbers.
Let $X$ be a compact K\"ahler manifold and
$g \in \Aut(X)$ an automorphism.
The pair $(X, g)$ is {\it strongly primitive}
if it is not bimeromorphic to another pair $(Y, g_Y)$
(even after replacing $g$ by its power)
having an equivariant fibration $Y \to Z$ with $\dim Y > \dim Z > 0$.
$g$ is of {\it positive entropy} if its {\it topological entropy} 
$$h(g) := \, \max \, \{ \log |\lambda| \,\,; \,\,
\lambda \,\,\,\, \text{is an eigenvalue of} \,\,\,\, g^* | \oplus_{i \ge 0} H^i(X, \BCC) \}$$
is positive; see \ref{conv}.
We remark that every surface automorphism of positive entropy
is automatically strongly primitive (cf. Lemma \ref{irrat}).

Theorems \ref{ThA}, \ref{ThB} and \ref{ThD} are our
main results, where the latter determines the geometrical structure
for those compact K\"ahler $X$ with a strongly primitive automorphism.
A subvariety $B \subset X$ is {\it $g$-periodic} if $g^s(B) = B$ for some $s > 0$.
Let $\rho(X)$ be the {\it Picard number} of $X$.

\begin{theorem}\label{ThA}
Let $X$ be a compact K\"ahler manifold, and $g \in \Aut(X)$
a strongly primitive automorphism of positive entropy.
Then we have:
\begin{itemize}
\item[(1)]
$X$ has at most $\rho(X)$ of $g$-periodic prime divisors.
\item[(2)]
If $X$ is a smooth projective threefold, then
any prime divisor of $X$ containing infinitely many
$g$-periodic curves, is itself $g$-periodic.
\end{itemize}
\end{theorem}

\begin{remark}\label{rThA}
\begin{itemize}
\item[(1)]
Suppose that the $X$ in Theorem \ref{ThA} (1) has $\rho(X)$ of $g$-periodic prime divisors,
then the algebraic dimension $a(X) = 0$
by the proof, Theorem \ref{ThD} and Remark \ref{rstabD}. Suppose further that
the irregularity
$q(X) := h^1(X, \SO_X) > 0$. Then the albanese map $\alb_X : X \to \Alb(X) = : Y$ is surjective and isomorphic
outside a few points of $Y$, and $\rho(Y) = 0$. Conversely,
we might realize such maximal situation by taking a complex $n$-torus $T$
with $\rho(T) = 0$
and a matrix $H \in \SL_n(\BZZ)$ with trace $> n$
so that $H$ induces an automorphism $h \in \Aut(T)$ of positive entropy;
if $H$ could be so chosen that $h$ has a few finite orbits $O_i$ of a total $\rho$ points $P_{ij} \in T$,
then the blowup $a: X \to T$ along these $\rho$ points lifts $h$ to some $g \in \Aut(X)$
of positive entropy with $\rho = \rho(X)$ of $g$-periodic prime divisors $a^{-1}(P_{ij})$.

\item[(2)]
When $\dim X = 2$, see \cite[Proposition 3.1]{Kawaguchi} or \cite[Theorem 6.2]{Z1} for
results similar to Theorem \ref{ThA} (1). Meromorphic endomorphisms and fibrations are studied in
\cite{AC}.

\item[(3)]
For a possible generalization of Theorem \ref{ThA}
to varieties over other fields,
we remark that the Bertini type theorem is
used in the proof, so the ground field might need to be
of characteristic zero.
K\"ahler classes are also used in the proof.
The proof of Theorem \ref{ThA} (2) requires
$X$ to be projective in order to define nef reduction as in \cite{8aut}.
\end{itemize}
\end{remark}

The following consequence of Lemma \ref{torus} or Theorem \ref{ThD}
and Lefschetz's fixed point formula,
shows the practicality of the strong primitivity notion.

\begin{theorem}
Let $A$ be a complex torus of $\dim A \ge 2$,
and $g \in \Aut_{\variety}(A)$ a strongly primitive automorphism of positive
entropy {\rm (cf. \ref{conv})}. Then $A$ has no $g$-periodic subvariety $D$ with
\ $\rm{pt} \ne D \subset A$.
In particular, for every $s > 0$, the number $\#\Per(g^s)$ of
$g^s$-fixed points $($with multiplicity counted$)$ satisfies
$$\#\Per(g^s) \, = \, \sum_{i \ge 0} \,\, \Tr \, (g^s)^* \, | \, H^i(A, \BZZ).$$
\end{theorem}

\par \noindent
{\bf Acknowledgement}

I thank Tien-Cuong Dinh and Nessim Sibony for the informative reference \cite{DS08} and comments
about $g$-periodic points, and the referee for the better exposition of the paper.

\section{Preliminary results}

\begin{setup}\label{conv}
{\bf Most of the conventions} are as in \cite{KM} and Hartshorne's book. Below are some more.

In the following (till Lemma \ref{irrat}), $X$ is a compact
K\"ahler manifold of dimension $n \ge 2$.

(1) Denote by $\NS(X) = \Pic(X)/\Pic^0(X)$ the
{\it Neron-Severi group}, and $\NS_{B}(X) = \NS(X) \otimes_{\BZZ} B$
for $B = \BQQ, \BRR$,
which is a $B$-vector space of finite dimension $\rho(X)$ (called the {\it Picard number}).

By abuse of notation, the {\it cup product} $L \cup M$ for
$L \in H^{i,i}(X)$ and $M \in H^{j,j}(X)$ will be denoted as
$L . M$ or simply $L M$.
Two codimension-$r$ cycles $C_1, C_2$ are {\it numerically equivalent} if
$(C_1 - C_2) M_1 \cdots M_{n-r}$ $= 0$ for all $M_i \in H^{1,1}(X)$.
Denote by $[C_1]$ the equivalence class containing $C_1$,
and $N^{r}(X)$ the $\BRR$-vector space of all equivalence classes $[C]$
of codimension $r$-cycles.
By {\it abuse of notation}, we write $C_1 \in N^r(X)$ (instead of $[C_1] \in N^r(X)$).
We remark that if $C_1$ and $C_2$ are cohomologous then
$C_1$ and $C_2$ are numerically equivalent, but the converse may not be true
if $r \le n-2$. Our $N^{n-1}(X)$ coincides with the usual $N_1(X)$.

Codimension-$r_i$ cycles $C_i$ ($i = 1, 2$) are {\it perpendicular} to each other
if $C_1 . C_2 = 0$ in $N^{r_1+r_2}(X)$.

(2) A class $L$ in the closure of the
K\"ahler cone of $X$ is called {\it nef}; this $L$ is {\it big} if $L^{n} \ne 0$.

For $g \in \Aut(X)$, the {\it $i$-th dynamical degree} is defined
as 
$$d_i(g) := \max \, \{|\lambda| \, ; \, \lambda \,\,\,\, \text{is
an eigenvalue of} \,\,\,\,  g^* | H^{i, i}(X) \} .$$
It is known that the {\it topological entropy} $h(g)$ equals $\max_{1 \le i \le n} \log d_i(g)$.
We say that $g$ is of {\it positive entropy} if $h(g) > 0$.
Note that $h(g) > 0$ if and only if $d_i(g) > 1$
for some $i$ and in fact for all $i \in \{1, \dots, n-1\}$, if and only if $h(g^{-1}) > 0$.
We refer to \cite{DS} for more details.

By the generalized Perron-Frobenius theorem in \cite{Bi}, there are nonzero nef classes
$L_g^{\pm}$ such that $g^*L_g^+ = d_1(g) L_g^+$ and $(g^{-1})^*L_g^- = d_1(g^{-1}) L_g^-$
in $H^{1,1}(X)$.
When $X$ is a projective manifold, we can choose $L_g^{\pm}$ to be in
$\NS_{\BRR}(X)$.

An irreducible subvariety $Z$ of $X$ is $g$-{\it periodic} if $g^s(Z) = Z$ for some $s \ge 1$.

(3) When a cyclic group $\langle g \rangle$ acts on $X$,
we use $g|X$ or $g_X$ to denote the image of $g$ in $\Aut(X)$.
The pair $(X, g|X)$ is loosely denoted as $(X, g)$.

(4) Suppose that a cyclic group $\langle g \rangle$ acts on compact
K\"ahler manifolds $X$, $X_i$, $Y_j$.
A morphism $\sigma : X_1 \to X_2$ is $g$-{\it equivariant}
if $\sigma \circ g = g \circ \sigma$.
Two pairs $(Y_1, g)$ and $(Y_2, g)$ are {\it bimeromorphically equivariant}
if there is a decomposition $Y_1 = Z_1 \overset{\sigma_1}{\ratmap} Z_2 \cdots
\overset{\sigma_r}{\ratmap} Z_{r+1} = Y_2$ into bimeromorphic maps such that for each $i$ either
$\sigma_i$ or $\sigma_i^{-1}$ is a $g$-equivariant bimeromorphic morphism.

$(X, g)$ or simply $g|X$, is {\it non-strongly-primitive} (resp. {\it non-weakly-primitive})
if $(X, g^s)$, for some $s > 0$,
is bimeromorphically equivariant to some
$(X', g^s)$ and there is a $g^s$-equivariant surjective morphism $X' \to Z$
with $Z$ a compact K\"ahler manifold of $\dim X > \dim Z > 0$
(resp. of $\dim X > \dim Z > 0$ and $g^s|Z = \id$).
We call $(X, g)$ {\it strongly primitive} (resp. {\it weakly primitive})
if $(X, g)$ is not non-strongly-primitive (resp. not non-weakly-primitive).

(5) For a complex torus $A$, the (variety) automorphism group $\Aut_{\variety}(A)$ equals
$T_A \rtimes \Aut_{\group}(A)$,
with $T_A$ the group of translations and $\Aut_{\group}(A)$ the group
of group-automorphisms.
\end{setup}

We frequently use the (5) below.
In particular, bimeromorphically equivariant automorphisms have
the same dynamical degrees (and hence entropy).

\begin{lemma}\label{irrat}
Let $X$ be a compact K\"ahler manifold of dimension $n$,
and $g \in \Aut(X)$ an automorphism of positive entropy.
Then the following are true.
\begin{itemize}
\item[(1)]
We have $n \ge 2$. If $n = 2$, then
$g$ is strongly primitive.
\item[(2)]
All $d_i(g^{\pm})$ $(1 \le i \le n-1)$ are irrational algebraic integers.
\item[(3)]
Let $L_i$ $(1 \le i \le n-1)$ be in the closure $\overline{P^i(X)}$ of the K\"ahler cone 
$P^i(X)$ of degree $i$
in the sense of \cite[Appendix A, Lemma A.9, the definition before Lemma A.3]{NZ}
such that $g^*L_i = d_i(g) L_i$ in $H^{i,i}(X)$. Then
no positive multiple of $L_i$ is in $H^{2i}(X, \BQQ)$.
\item[(4)]
Every $g$-periodic curve is perpendicular to $L_1$.
\item[(5)]
We have $d_i(g) = d_i(g|Y)$ $(1 \le i \le n)$ if there is a $g$-equivariant
generically finite surjective morphism either from $X$ to $Y$ or
from $Y$ to $X$.
Here $g$ is not assumed to be of positive entropy.
\end{itemize}
\end{lemma}

\begin{proof}
For (1), apply Lemma \ref{DSl} or \cite[Lemma 2.12]{Z1}
to $L_g^+$ and the fibre of an equivariant fibration (cf. also (5)).
For the existence of the $L_i$ in (3), we used the generalized Perron-Frobenius theorem 
in \cite{Bi} for the closed cone $\overline{P^i(X)} \subset H^{i, i}(X, \BRR)$.
Now (3) follows from (2) by considering the cup product.

(2) Since $g^{-1}$ is also of positive entropy, we consider only $g$.
Since $g^*$ acts on $H^i(X, \BZZ)$ and each $d_i(g)$ is known to be an eigenvalue
of $H^i(X, \BCC) = H^i(X, \BZZ) \otimes_{\BZZ} \BCC$, all
dynamical degrees $d_i(g) > 1$ are algebraic integers.
Suppose that $d_i(g)$ is rational. Then $d_i(g) \in \BZZ_{\ge 2}$.
Take an eigenvector $M_i$ in $H^{2i}(X, \BZZ)$ with $g^*M_i = d_i(g)M_i$.
Since the cup product is non-degenerate, we can find $N_{n-i} \in H^{2n-2i}(X, \BZZ)$
such that $M_i . N_{n-i} = m_i \in \BZZ \setminus \{0\}$.
Now $m_i/d_i(g)^s = (g^{-s})^*M_i . N_{n-i} \in \BZZ$
for all $s > 0$. This is absurd.

(4) Suppose that $g^s(C) = C$ for some $s > 0$ and a curve $C$. Then
$L_1 . C = (g^s)^*L_1 . (g^s)^*C = d_1(g)^s L_1 . C$.
So $L_1 . C = 0$ for $d_1(g) > 1$.

For (5), see \cite[Lemma 2.6]{Z2} and \cite[Appendix A, Lemma A.8]{NZ}.
\end{proof}

The result below should be well known (cf.~e.g.~\cite[Appendix, Lemma A.4]{NZ}).

\begin{lemma}\label{comp0}
Let $X$ be a compact K\"ahler manifold and $L \in \overline{P^i(X)}$ 
{\rm (cf.~Lemma \ref{irrat} (3) for the notation)}.
Then $L = 0$ in $N^i(X)$ if and only if $L = 0$ in $H^{i,i}(X, \BRR)$.
\end{lemma}

The two results below are crucial and due to Dinh-Sibony \cite{DS}, but we slightly reformulated
(cf.~Lemma \ref{comp0}).

\begin{lemma} {\rm(cf.~\cite[Lemme 4.4]{DS})} \label{DSl}
Let $X$ be a compact K\"ahler manifold of dimension $n \ge 2$,
$g : X \to X$ a surjective endomorphism, and $M_1, M_2, L_i$ $(1 \le i \le m; m \le n-2)$
nef classes. Suppose that in $N^{m+1}(X)$ we have $L_1 \cdots L_m M_i \ne 0$ $(i = 1, 2)$
and $g^*(L_1 \cdots L_m M_i) = \lambda_i (L_1 \cdots L_m M_i)$ for some $($positive real$)$
constants $\lambda_1 \ne \lambda_2$. Then $L_1 \cdots L_m M_1 M_2 \ne 0$
in $N^{m+2}(X)$.
\end{lemma}

\begin{lemma}{\rm(cf.~\cite[Corollaire 3.2]{DS})} \label{DSc}
Let $X$ be a compact K\"ahler manifold with nef classes $L, M$.
If $L M = 0$ in $N^2(X)$, then $L$ and $M$ are parallel in $H^{1, 1}(X, \BRR)$.
\end{lemma}

\begin{lemma}\label{n-1}
Let $X$ and $Y$ be compact K\"ahler manifolds with $n:= \dim X \ge 2$,
and $\pi : (X, g) \to (Y, g_Y)$ an equivariant surjective morphism.
\begin{itemize}
\item[(1)]
Suppose that a nef and big class $M$ on $X$ satisfies $g^*M = M$ in $H^{1,1}(X)$.
Then a positive power of $g$ is in $\Aut_0(X)$ and hence
$g$ is of null entropy.
\item[(2)]
Suppose that $g$ is of positive entropy and $\dim Y = n-1$.
Then no nef and big class $M$ on $Y$ satisfies $g_Y^*M = M$.
In particular, $g_Y^* | H^{1,1}(Y)$ is of infinite order
and hence no positive power of $g_Y$ is in the identity
connected component $\Aut_0(Y)$ of $\Aut(Y)$.
\end{itemize}
\end{lemma}

\begin{proof}
(1) is a result of Lieberman \cite[Proposition 2.2]{Li}; see \cite[Lemma 2.23]{Z2}
(by Demailly-Paun, a nef and big class can be written as the sum
of a K\"ahler class and a closed real positive current).

(2) If $g_Y^* | H^{1,1}(Y)$ is of finite order $r$, then
$g_Y^*$ stabilizes $\sum_{i=0}^{r-1} (g_Y^i)^*H$ with $H$ a K\"ahler class.
So we only need to rule out the existence of such $M$ in the first assertion.
Set $M_X := \pi^*M$. We apply Lemma \ref{DSl} repeatedly to
show the assertion that $M_X^{k-1} . L_g^+ \ne 0$ in $N^k(X)$ for all $1 \le k \le n$.
Indeed, $M_X . L_g^+$ is nonzero in $N^2(X)$ since $g^*M_X = M_X$ while $g^*L_g^+ = d_1(g) L_g^+$
with $d_1(g) > 1$; if $M_X^{j-1} . L_g^+ \ne 0$ in $N^j(X)$ for $j < n$, then
$M_X^j . L_g^+ \ne 0$ in $N^{j+1}(X)$ because $g^*(M_X^{j-1} . L_g^+) = d_1(g)(M_X^{j-1} . L_g^+)$
with $d_1(g) > 1$, and $g^*M_X^j = M_X^j$ ($\ne 0$ in $N^j(X)$), so the assertion is true.
Now $\deg(g) (M_X^{n-1} . L_g^+) = g^*M_X^{n-1} . g^*L_g^+ = d_1(g)(M_X^{n-1} . L_g^+)$
implies a contradiction:
$1 = \deg(g) = d_1(g) > 1$.
Lemma \ref{n-1} is proved.
\end{proof}

\begin{lemma}\label{LgNullD}
Let $X$ be a compact K\"ahler manifold of dimension $n \ge 2$ and $q(X) = 0$, and $g \in \Aut(X)$
an automorphism of positive entropy.
Then $X$ has at most $\rho(X)$ of prime divisors $D_j$
perpendicular to either one of $L_g^+$ and $L_g^-$
in $N^2(X)$. Further, such $D_j$ are all $g$-periodic.
\end{lemma}

\begin{proof}
We only need to show the first assertion, since
both $L_g^{\pm}$ are semi $g^*$-invariant and hence $g$ permutes these $D_j$.

Suppose that
$X$ has $1 + \rho(X)$ of distinct prime divisors
$D_i$ with $L_g^+ . D_i = 0$
in $N^2(X)$.
The case $L_g^-$ is similar by considering $g^{-1}$. Set
$L := L_g^+$.
Since these $D_i$ are then linearly dependent, we may assume that
$E_1 := \sum_{i=1}^{t_1} a_i D_i \equiv E_2 := \sum_{j=t_1+1}^{t_1+t_2} b_j D_j$
in $\NS_{\BQQ}(X)$
for some positive integers $a_i, b_j, t_k$.
Since $q(X) = 0$, we may assume that $E_1 \sim E_2$ (linear equivalence)
after replacing $E_i$ by its multiple.
Let $\sigma: X' \to X$ be a blowup such that $|\sigma^*E_1| = |M| + F$
with $|M|$ base point free and $F$ the fixed component.
Take a K\"ahler class $H$ on $X$. Then $0 \le \sigma^*L . M . \sigma^*(H^{n-2}) \le
\sigma^*L . (M + F) . \sigma^*(H^{n-2}) = L . E_1 . H^{n-2} = 0$.
Hence $\sigma^*L . M . \sigma^*(H^{n-2}) = 0$. Thus, $\sigma^*L . M = 0$
in $H^{2,2}(X', \BRR)$ by \cite[Appendix A, Lemmas A.4 and A.5]{NZ}.
So, by Lemma \ref{DSc}, $\sigma^*L$ equals $M$ in $\NS_{\BQQ}(X')$,
after replacing $L$ by its multiple.
Thus $L \in \NS_{\BQQ}(X)$, contradicting
Lemma \ref{irrat}.
This proves Lemma \ref{LgNullD}.
\end{proof}

Theorem \ref{stabD} below effectively bounds the number of $g$-periodic prime divisors.

\begin{theorem}\label{stabD}
Let $X$ be a
compact K\"ahler
manifold of dimension $n \ge 2$ and $q(X) = 0$, and $g \in \Aut(X)$
a weakly primitive automorphism of positive entropy.
Then we have:
\begin{itemize}
\item[(1)]
$X$ has none or only finitely many $g$-periodic prime divisors
$D_i$ $(1 \le i \le  r; \, r \ge 0)$.
\item[(2)]
If $r > \rho(X)$, then $n \ge 3$ and $($after replacing $g$ by its power
and $X$ by its $g$-equivariant blowup$)$ there is an equivariant surjective
morphism $\pi: (X, g) \to (Y, g_Y)$ with connected fibres,
$Y$ rational and almost homogeneous,
$\dim Y \in \{1, \dots, \, n-2\}$, and $g_Y \in \Aut_0(Y)$.
\item[(3)]
If $g$ is strongly primitive, then
$X$ has at most $\rho(X)$ of $g$-periodic prime divisors.
\end{itemize}
\end{theorem}

\begin{proof}
Let $D_i$ ($1 \le i \le r; \, r > \rho := \rho(X)$)
be distinct $g$-periodic prime divisors of $X$.
Then $D_i$'s are linearly dependent.
Replacing $g$ by its power, we may assume that $g(D_i) = D_i$ for all $i \le r$.
By the reasoning in Lemma \ref{LgNullD}, the Iitaka $D$-dimension
$\kappa := \kappa(X, \sum_{i = 1}^r D_{i}) \ge 1$.
If $\kappa = n$, then replacing $X$ by its $g$-equivariant blowup,
we may assume that some positive combination $M$ of $D_i$ is nef and big
and $g^*M = M$, contradicting Lemma \ref{n-1}.
Thus, $1 \le \kappa < n$.

Take $E_1 := \sum_{i=1}^t a_i D_i$ with $a_i$ non-negative integers
such that $\Phi_{|E_1|} : X \ratmap \BPP^N$ has the image $Y$ with
$\dim Y = \kappa$, and the induced map $\pi : X \ratmap Y$ has connected general fibres.
Since $g(E_1) = E_1$,
replacing $X$ by its $g$-equivariant blowup
and removing redundant components in $E_1$,
we may assume that $\Bs|E_1| = \emptyset$, $\pi$ is holomorphic, $Y$ is smooth projective,
and $g$ descends to an automorphism $g_Y \in \Aut(Y)$;
further we can write $E_1 = \pi^*A$, where $g_Y(A)$ equals $A$ and is a nef and
big Cartier divisor with $\Bs|A| = \emptyset$
(notice that $A$ may not be ample because we have replaced $Y$ by its blowup).
Hence
$g_Y \in \Aut_0(Y)$ after $g$ is replaced by its power, so
$\dim Y \ne n-1$; see Lemma \ref{n-1}.
Therefore, $1 \le \kappa = \dim Y \in \{1, \dots, n-2\}$.

By the assumption on $g$, we have $\ord(g_Y) = \infty$.
Since $q(Y) \le q(X) = 0$, our $\Aut_0(Y)$ is a linear algebraic group;
see \cite[Theorem 3.12]{Li} or \cite[Corollary 5.8]{Fu}.
Let $H$ be the identity component of the closure of $\langle g_Y \rangle$ in $\Aut_0(X)$,
and we may assume that $g_Y \in H$ after replacing $g$ by its power.
Let $\tau: Y \ratmap Z = Y/H$ be the quotient map; see \cite[Theorem 4.1]{Fu}.
Replacing $Y, Z, X$ by their equivariant blowups, we may assume that
$Y$ and $Z$ are smooth and $\tau$ is holomorphic.
By the construction,
$g \in \Aut(X)$ and $g_Y \in \Aut(Y)$ descend to $\id_Z \in \Aut(Z)$.
The assumption on $g$ implies that $\dim Z = 0$. So $Y$ has
a Zariski-open dense $H$-orbit $Hy$. In other words, $Y$ is almost homogeneous.
Since $H$ is abelian (and a rational variety by a result of Chevalley),
$Y$ is bimeromorphically dominated by $H$
(each stabilizer subgroup $H_y$ being normal in $H$),
so $Y$ is rational (and smooth projective).
(2) and (3) are proved.

To prove (1),
suppose that $X$ has infinitely many distinct $g$-periodic prime divisors
$D_i$ ($i \ge 1$). We may assume that $\kappa := \kappa(X, \sum_{i=1}^r D_i) =
\max\{\kappa(X, \sum_{i=1}^s D_i) \, | \, s \ge 1\} \ge 1$ for some $r > 0$,
and use the notation above. In particular, $1 \le \kappa \le n-2$.
We assert that $(*)$ all $D_j$ ($j > r$) are mapped to
distinct $g_Y$-periodic prime divisors $D_j' \subset Y$
by the map $\pi : X \to Y$, after replacing $\{D_i\}$ by an infinite subsequence.
Since $\pi$ is smooth (and hence flat) outside a codimension one subset of $X$
and the $\pi$-pullback of a prime divisor has only finitely many irreducible components,
we have only to consider the case where $D_{j_1}, D_{j_2}, \dots$ (with $j_v > r$)
is an infinite sequence of divisors each dominating $Y$, and show that this case is impossible.
Replacing $g$ by its power and $X$ by its $g$-equivariant blowup,
we may assume that $|E_3|$ is base point free for some
$E_3 = b_{j_1}D_{j_1} + \dots + b_{j_u} D_{j_u}$ with $b_{j_v} \in {\BZZ}_{\ge 1}$,
and $D_{j_1}$ dominates $Y$ (notice that some components of $E_3$ are in the
exceptional locus of the blowup).
By the maximality of $\kappa$, we have $\kappa(X, E_1 + E_3) = \kappa(X, E_1)$
and hence $\Phi_{|E_1+ E_3|}$ is holomorphic onto a variety $W$ of dimension $\kappa$
with $E_1 + E_3$ the pullback of an ample divisor $A_W \subset W$.
Thus taking a K\"ahler class $M$ on $X$, we obtain
a contradiction:
$$
0 = M^{n-1-\kappa} (E_1 + E_3)^{\kappa + 1} \ge M^{n-1-\kappa} . E_1^{\kappa} . E_3 \ge
M^{n-1-\kappa} . E_1^{\kappa} . D_{j_1} \\
= M^{n-1-\kappa} . B \, > \, 0
$$
where $E_1 = \pi^*A$ with $A$ nef and big as above, and $B = (\pi^*A | D_{j_1})^{\kappa}$
is a sum of $A^k$ of $(n-1-\kappa)$-dimensional general fibres of 
the surjective morphism $\pi|D_{j_1} : D_{j_1} \to Y$.
The assertion $(*)$ is proved.

Now the infinitely many distinct $g_Y$-periodic prime divisors $D_j' \subset Y$
are squeezed in
the complement of some Zariski-open dense $H$-orbit $Hy$ of $Y$ (for some general $y \in Y$,
whose existence was mentioned early on). 
This is impossible.
Thus, we have proved (1). The proof of
Theorem \ref{stabD} is completed.
\end{proof}

\begin{remark} \label{rstabD}
Assume that the algebraic dimension $a(X) = \dim X$ in Theorem \ref{stabD}. 
Then $X$ is projective since $X$
is K\"ahler. If $X$ has $\rho(X)$ of linearly independent
$g$-periodic divisors, then (a power of) $g^*$ stabilizes an ample divisor on
$X$; so $g$ is of null entropy by Lemma \ref{n-1}, absurd!
Thus, by the proof, `$r > \rho(X)$' in Theorem \ref{stabD} (2)
(resp. `$\rho(X)$' in Theorem \ref{stabD} (3)) can be replaced by
`$r \ge \rho(X)$' (resp. `$\rho(X) - 1$').
\end{remark}

\begin{lemma}\label{nefR}
Let $X$ be a projective manifold of dimension $n \ge 2$, and $g \in \Aut(X)$
an automorphism of positive entropy. Let $L = L_g^+$ or $L_g^-$.
Then the nef dimension $n(L) \ge 2$, and
the nef reduction map $\pi : X \ratmap Y$ in \cite{8aut} can be taken to be
holomorphic with $Y$ a projective manifold, after $X$ is replaced by its
$g$-equivariant blowup.
\end{lemma}

\begin{proof}
Since $L \ne 0$, we have $n(L) = \dim Y \ge 1$.
The second assertion is true by the construction of the nef reduction in \cite[Theorem 2.6]{8aut},
using the chain-connectedness equivalence relation
defined by numerically $L$-trivial curves (and preserved by $g$).
Consider the case $n(L) = 1$. For a general fibre $F$ of $\pi$,
we have $L|F = 0$ by the definition of the nef reduction. By Lemma \ref{DSc},
a multiple of $L$ is equal to $F$ in $\NS_{\BQQ}(X)$, contradicting Lemma \ref{irrat}.
\end{proof}

We remark that the hypothesis in Lemma \ref{LggNull} below is optimal and
the hypothetical situation may well occur when $X \to Y$ is $g$-equivariant,
$Y$ is a surface, and $D_j$ and $L_g^{\pm}$ are
pullbacks from $Y$, e.g. when $X = Y \times$ (a curve) and $g = g_Y \times \id$.

\begin{lemma}\label{LggNull}
Let $X$ be a $3$-dimensional projective manifold with $q(X) = 0$, and $g \in \Aut(X)$
an automorphism of positive entropy. Let $D_i$ $(i \ge 1)$ be infinitely many
pairwise distinct prime divisors such that
$L_g^+ . L_g^- . D_i = 0$. Then for both $L = L_g^+$ and $L = L_g^-$, we have
$L^2 = 0$ in $N^2(X)$ and the nef dimension
$n(L) = 2$.
\end{lemma}

\begin{proof}
Note that $L_g^+ . L_g^- \ne 0$ in $N^2(X)$ by Lemma \ref{DSl} or \ref{DSc}.
Set $L_1 := L_g^+$, $L_2 := L_g^-$ and
$\lambda_1 := d_1(g) > 1$, $\lambda_2 := 1/d_1(g^{-1}) < 1$.
Then $g^*L_i = \lambda_i L_i$.
If $L_i^2 \ne 0$ in $N^2(X)$ for both $i = 1, 2$,
then $L_i . L_i . L_j \ne 0$, where
$\{i, j\} = \{1, 2\}$; see Lemma \ref{DSl};  applying $g^*$, we get
$\lambda_i^2 \lambda_j = 1$, whence $1 < \lambda_1 = \lambda_2 < 1$, absurd.

To finish the proof of the first assertion,
we only need to consider the case where $L_1^2 \ne 0$ and $L_2^2 = 0$
in $N^2(X)$, because we can switch $g$ with $g^{-1}$.
By Lemma \ref{DSl}, $L_1^2 . L_2 \ne 0$.
Now $L_1 + L_2$ is nef and big because
$(L_1+L_2)^3 \ge 3L_1^2 L_2 > 0$.
So we can write $L_1 + L_2 = A + \Delta$ with an ample $\BRR$-divisor $A$ and an effective
$\BRR$-divisor $\Delta$; see \cite[Lemma 2.23]{Z2}
for the reference on such decomposition.
By Lemma \ref{LgNullD} and taking an infinite subsequence,
we may assume that
$L_i . D_j \ne 0$ in $N^2(X)$ for $i = 1$ and $2$ and all $j \ge 1$, and
$D_j$ is not contained in the support of $\Delta$ for all $j \ge 1$.
Now $L_1^2 . D_j = (L_1 + L_2)^2 . D_j = (L_1 + L_2) . (A + \Delta) . D_j
\ge (L_1 + L_2) . A . D_j \ge A^2 . D_j > 0$. Thus $L_1|D_j$ is a nef and big divisor and
$L_2|D_j$ is a nonzero nef divisor such that $(L_1|D_j) . (L_2|D_j) = L_1 . L_2 . D_j = 0$.
This contradicts the Hodge index theorem applied to a resolution of $D_j$.
The first assertion is proved.

Let $L$ be one of $L_g^+$ and $L_g^-$. By Lemma \ref{nefR}, we only need to show $n(L) \ne 3$.
As in the proof of Theorem \ref{stabD}, we may assume that the Iitaka $D$-dimension
$\kappa := \kappa(X, E_1) = \max\{\kappa(X, \sum_{i=1}^s D_{i}) \, | \, s \ge 1\} \ge 1$
for some $E_1 := \sum_{i=1}^{t} a_i D_i$ with positive integers $a_i$.
If $\kappa(X, E_1) = 3$, then $E_1$ is big and hence a sum of an ample
divisor and an effective divisor, whence $L_g^+ . L_g^- . E_1 > 0$,
contradicting the choice of $D_j$. Therefore,
$\kappa = 1, 2$.

Case (1). $\kappa = 2$. Let $\sigma : X' \to X$ be a blowup such that
$|\sigma^*E_1| = |M| + F$ with $|M|$ base point free and $F$ the fixed component.
Since $\kappa(X', M) = \kappa(X, E_1) = 2$, we have $M^2 \ne 0$.
If $\sigma^*L . M^2 = 0$, then the projection formula implies that
$L . C = 0$ for every curve $C = \sigma_*(M_1 . M_2)$
with $M_i \in |M|$ general members. So the nef dimension $n(L) < 3$.

Suppose that $\sigma^*L . M^2 > 0$. Then $\sigma^*L + M$
is nef and big because $(\sigma^*L+M)^3 \ge 3\sigma^*L . M^2 > 0$.
Since $\sigma^*(L + E_1)$ is larger than $\sigma^*L + M$,
it is also big. So $L + E_1$ is big, too.
Hence $0 < L . L' . (L + E_1) = L_g^+ . L_g^- . E_1$,
where $\{L, L'\} = \{L_g^{\pm}\}$,
contradicting
the choice of $D_j$ and $E_1$.

Case (2). $\kappa = 1$. We may assume that $|E_1|$ has no fixed
component and is an {\it irreducible}
pencil parametrized by $\BPP^1$ (noting: $q(X) = 0$), after removing redundant $D_j$ from $E_1$.
Since $L_g^{\pm}$ are semi $g^*$-invariant, every $g(D_j)$, like $D_j$, is
also perpendicular to $L_g^+ . L_g^-$.
After relabelling and expanding the sequence, we may assume that
$g(E_1)$ is also a positive combination of $D_j$'s.
By Case (1), we may assume that $\kappa(E_1 + g(E_1)) = 1$.
For general (irreducible) members $M_1 \in |E_1|$ and
$M_2 \in |g(E_1)|$, the two-component divisor $M_1 + M_2$ is a reduced member of
$|E_1 + g(E_1)|$.

Note that $N:= h^0(E_1 + g(E_1)) \ge h^0(E_1) + h^0(g(E_1)) - 1 \ge 3$.
The linear system $|E_1 + g(E_1)|$ gives rise to a rational map from $X$
onto a curve $B$ of degree $\ge N-1$ in $\BPP^{N-1}$.
Thus, each member of $|E_1 + g(E_1)|$ lying over $B \setminus \Sing B$,
is a sum of $N-1$ linearly equivalent nonzero effective divisors, since $B$ is a rational
curve; indeed, the genus $g(B)$ of $B$ satisfies $g(B) \le q(X) = 0$. So $E_1 \sim g(E_1)$.
Replacing $X$ by its $g$-equivariant blowup,
we may assume that $|E_1|$ is base point free and hence $E_1$ is a nef eigenvector of $g^*$.
Now $L_g^+ . L_g^- . E_1 = 0$ infers a contradiction to
Lemma \ref{DSl}, since $L_g^+$, $L_g^-$ and $E_1$ correspond to distinct eigenvalues $d_1(g), 1/d_1(g^{-1}), 1$ of $g^*|\NS_{\BQQ}(X)$.
This proves Lemma \ref{LggNull}.
\end{proof}

\begin{lemma}\label{torus}
Let $A$ be a complex torus of dimension $n \ge 2$
and $f \in \Aut_{\variety}(A)$ of infinite order such that
$f(D) = D$ for some subvariety \ $\rm{pt} \ne D \subset X$.
Then there is a subtorus $B \subset A$ with $\dim B \in \{1, \dots, n-1\}$
such that $f$
descends, via the quotient map $A \to A/B$,
to an automorphism $h \in \Aut_{\variety}(A/B)$
having a periodic point in $A/B$.
\end{lemma}

\begin{proof}
Write $f = T_a \circ g$ with $T_a \in T_A$ a translation and $g$
a group automorphism.

Case (1). $\kappa(D) = \dim D$, i.e.,
$D$ is of general type. Then $\Aut(D)$ is finite, so
$f^s|D = \id_D$ for some $s > 0$.
Since $f^s$ fixes
$D$ pointwise,
the identity component $B$ of the pointwise fixed point set
$A^{g^s}$ (a subtorus)
is a positive-dimensional subtorus; see \cite[Lemma 13.1.1]{BL}.
Write $f^s = T_c \circ g^s$ with $T_c \in T_A$.
If $\dim B \ge n$, then $B = A$, $g^s = \id_A$ and $f^s = T_c$, so $f^s = \id$
for $f^s|D = \id_D$. This contradicts
the assumption on $f$. Thus $1 \le \dim B \le n-1$.
Our $g$ acts on $A^{g^s}$, so $g(B) \subset A^{g^s}$ is a coset
in $A^{g^2}/B \le A/B$. Thus $g(B) = \delta + B$ for some $\delta$.
So $g(B) = B$, because $(*)$ : $g$ is a group-automorphism
and $0 \in B \le A$.
Now $f(x + B) = a + g(x) + g(B) = f(x) + B$.
So $f$ permutes cosets in $A/B$
and $f^s$ fixes those cosets $d + B$ with $d \in D$.
Lemma \ref{torus} is true.

Case (2). The Kodaira dimension $\kappa(D) \le 0$.
Then $\kappa(D) = 0$ and
$D = \delta + B$ with a subtorus $B$ of $A$;
see \cite[Lemma 10.1, Theorem 10.3]{Ue}.
Now $\delta + B = D = f(D) = a + g(\delta) + g(B)$,
thus $g(B)$ equals a coset in $A/B$ and hence $g(B) = B$
by the reasoning (*) in Case (1).
Therefore, $f$ permutes cosets in $A/B$ as in Case (1),
and fixes the coset $\delta + B$. So
Lemma \ref{torus} is true.

Case (3). $\kappa(D) \in \{1, \dots, \, \dim D -1\}$.
By \cite[Theorem 10.9]{Ue},
the identity connected component $B$ of $B':= \{x \in A \, | \, x + D \subseteq D\}$
is a subtorus with $\dim B = \dim D - \kappa(D)$.
We claim that $f$ permutes cosets in $A/B$.
Indeed, for every $b \in B$, we have
$D = f(D) = f(b+D) = a + g(b) + g(D) = g(b) + f(D) = g(b) + D$,
so $g(b) \in B'$. Thus $g(B) \le B'$. Hence $g(B) = B$ and the claim is true,
by the reasoning in Case (1).
Further, the map $D \to D/B$ is bimeromorphic to the Iitaka fibration,
and $\kappa(D/B) = \dim (D/B)$ (cf. ibid.).
$f$ descends to an automorphism $f' \in \Aut_{\variety}(A/B)$
stabilizing $D/B \subset A/B$. Using Case (1), we are done
for some quotient torus $(A/B)/(B'/B) \cong A/B'$.
Lemma \ref{torus} is proved.
\end{proof}

\section{Proof of Theorem \ref{ThA}}

In this section, we prove Theorem \ref{ThA} in the introduction
and the two results below. Theorem \ref{ThB} treats $X$ with $q(X) = 0$,
while Theorem \ref{ThD} determines the geometrical structure of those K\"ahler $X$ with
a strongly primitive automorphism.

\begin{theorem}\label{ThB}
Let $X$ be a compact K\"ahler manifold of dimension $n \ge 2$ and irregularity
$q(X) = 0$, and $g \in \Aut(X)$
a weakly primitive automorphism of positive entropy. Then:
\begin{itemize}
\item[(1)]
$X$ has finitely many prime divisors $B_i$ $(1 \le i \le r; \, r \ge 0)$ such that:
each $B_i$ is $g$-periodic, and
$\cup B_i$ contains every $g$-periodic prime divisor and every prime divisor
perpendicular to $L_g^{+}$ or $L_g^-$.
\item[(2)]
Suppose that $g$ is strongly primitive.
Then the $r$ in $(1)$ satisfies $r \le \rho(X)$,
and $r = \rho(X)$ holds only when the algebraic dimension $a(X) < n$.
\item[(3)]
Suppose that $X$ is a smooth projective threefold, and $g$ is strongly primitive.
Then $(L_g^+ + L_g^-) | D$ is nef and big for every prime divisor $D \ne B_i$ $(1 \le i \le r)$.
In particular, if a prime divisor $D \subset X$ contains infinitely many curves
each of which is either $g$-periodic or perpendicular to $L_g^+ + L_g^-$,
then $D$ itself is $g$-periodic.
\end{itemize}
\end{theorem}

A compact K\"ahler manifold $X$ is called {\it weak Calabi-Yau}
if $\kappa(X) = 0 = q(X)$.

\begin{theorem}\label{ThD}
Let $X$ be a compact K\"ahler manifold of dimension $n \ge 2$, and $g \in \Aut(X)$
a strongly primitive automorphism of positive entropy.
Then the algebraic dimension $a(X) \in \{0, n\}$.
Suppose further that $(*)$
either $\kappa(X) \ge 0$, or $q(X) > 0$, or
$\kappa(X) = - \infty$, $q(X) = 0$ and $X$ is projective and uniruled.
Then $(1)$, $(2)$ or $(3)$ below occurs.
\begin{itemize}
\item[(1)]
$X$ is a weak Calabi-Yau manifold.
\item[(2)]
$X$ is rationally connected in the sense of
Campana, Kollar-Miyaoka-Mori $($so $q(X) = 0)$.
\item[(3)]
The albanese map $\alb_X : X \to \Alb(X)$ is surjective and
isomorphic outside a few points of $\Alb(X)$.
There is no $h$-periodic subvariety of dimension in $\{1, \dots, \, n-1\}$
for the $($variety$)$ automorphism $h$ of $\Alb(X)$ induced from $g$.
\end{itemize}
\end{theorem}

\begin{setup}
{\bf Proof of Theorem \ref{ThB}}
\end{setup}

The assertions (1) and (2) follow from Lemma \ref{LgNullD}, Theorem \ref{stabD}
and Remark \ref{rstabD}.
For (3), by Lemmas \ref{LggNull} and \ref{nefR}, our $X$ has finitely many divisors
$D_j$ ($1 \le j \le s$) such that $L_g^+ . L_g^- . D_j = 0$ and
$L_g^+ . L_g^- . D > 0$
for every prime divisor $D \ne D_j$ ($1 \le j \le s$).
Since both $L_g^{\pm}$ are semi $g^*$-invariant, these $D_j$'s are permuted
by $g$ and hence are all $g$-periodic.
Thus $\{D_j\} \subset \{B_i\}$.

Suppose that $D \ne B_i$ ($1 \le i \le r$) is a prime divisor of $X$.
Then $M:= L_g^+ + L_g^-$ is nef and
$(M|D)^2 \ge 2L_g^+ . L_g^- . D > 0$, so $M|D$ is nef and big.
Thus $D$ has none or only finitely many curves perpendicular to $M$,
by the Hodge index theorem applied to a resolution of $D$.
So $D$ contains only finitely many $g$-periodic curves (cf.~Lemma \ref{irrat} (4)).
This proves (3) and also Theorem \ref{ThB}.

\begin{setup}
{\bf Proof of Theorem \ref{ThD}}
\end{setup}

As in the proof of \cite[Lemma 2.16]{Z3}, a suitable algebraic reduction
$X \to Y$, with $\dim Y = a(X)$, is holomorphic and $g$-equivariant.
So $a(X) \in \{0, n\}$, since $g$ is strongly primitive.

Consider the case $\kappa(X) \ge 1$. Let
$\Phi = \Phi_{|mK_X|} : X \ratmap \BPP^{N}$ be the Iitaka fibration.
Replacing $X$ by its $g$-equivariant blowup, we may assume that
$\Phi$ is holomorphic and $g$-equivariant
onto some smooth $Z$ with $\dim Z = \kappa(X)$.
Our $g$ descends to an automorphism $g_Z \in \Aut(Z)$.
Now $\ord(g_Z) < \infty$ (so $\dim Z < \dim X$ by Lemma \ref{irrat} (5)), by the generalization of
\cite[Theorem 14.10]{Ue} to dominant
meromorphic selfmaps on K\"ahler manifolds as in \cite[Theorem A or Corollary 2.4]{NZ}.
This contradicts the strong primitivity of $g$.
Therefore, $\kappa(X) \le 0$.

Case(1). $q(X) > 0$. We will show that
Theorem \ref{ThD} (3) holds.
Consider the albanese map $\alb_X : X \to \Alb(X)$
and let $Y = \alb_X(X)$ be its image. $g$ descends to automorphisms $g | \Alb(X)$
and $h \in \Aut(Y)$. Since $g$ is strongly primitive,
$\dim Y = n$. Thus $\alb_X$ is generically finite onto $Y$
and hence $0 \ge \kappa(X) \ge \kappa(Y) \ge 0$; see \cite[Lemma 10.1]{Ue}.
So $\kappa(X) = \kappa(Y) = 0$.
Hence $\alb_X$ is surjective and bimeromorphic, with $E$ denoting the exceptional divisor;
see \cite[Theorem 24]{Ka}.
If $\alb_X$ is not an isomorphism, i.e., $E \ne \emptyset$, then $g(E) = E$ and
$h(\alb_X(E)) = \alb_X(E)$ because
$g$ and $h$ are compatible.
By Lemma \ref{torus} and since $g$ is strongly primitive, $\dim \alb_X(E) = 0$.
So Theorem \ref{ThD} (3) holds by Lemma \ref{torus}.

If $q(X) = 0 = \kappa(X)$, then $X$ is weak Calabi-Yau by the definition.
So we have only to consider the case where $q(X) = 0$ and $\kappa(X) = -\infty$,
or the following case by the assumption.

Case (2). $X$ is projective and uniruled. We will show that $X$ is rationally connected.
After $g$-equivariant blowups, we may assume that the
maximal rationally connected
fibration $\pi: X \to Y$ is holomorphic and $g$-equivariant, with $Y$ smooth
and $\dim Y < n$ (cf.~\cite[Theorem C]{NZ}).
Since $g$ is strongly primitive, we have $\dim Y = 0$,
so $X$ is rationally connected.
Theorem \ref{ThD} is proved.

\begin{setup}
{\bf Proof of Theorem \ref{ThA} and Remark \ref{rThA} (1)}
\end{setup}

For Theorem \ref{ThA} (1), by Theorem \ref{ThB}, we may assume that $q(X) > 0$,
so Theorem \ref{ThD} (3) occurs.
Suppose that $X$ has $r \ge \rho := \rho(X)$
of $g$-periodic prime divisors $D_i$.
Then each $\alb_X(D_{i}) \subset \Alb(X) = : Y$ is $h$-periodic,
so it is a point, since we are in
Theorem \ref{ThD} (3).
Thus these $D_i$ are irreducible components of
the exceptional divisor $E$ of $\alb_X : X \to Y$.
We assert that $(**) : \, \NS_{\BQQ}(X)$ has a basis consisting of
the irreducible components of $E$
and the pullback of a basis of $\NS_{\BQQ}(Y)$.
This is clear if $\alb_X$ is the blowup along a smooth centre.
The general case can be reduced to this special case by
the weak factorization theorem of bimeromorphic maps due to
Abramovich-Karu-Matsuki-Wlodarczyk (or by blowing up the indeterminacy
of $Y \ratmap X$ as suggested by Oguiso).
Now the assertion $(**)$ implies that $r = \rho$, $E = \sum_{i=1}^{\rho} D_i$
and $\rho(Y) = 0$ (so $a(X) = 0$ by Theorem \ref{ThD}).
This proves Theorem \ref{ThA} (1) and Remark \ref{rThA} (1).

For Theorem \ref{ThA} (2), let $D \subset X$ be a prime divisor containing
infinitely many $g$-periodic curves $C_i$ ($i \ge 1$).
We may assume that $q(X) > 0$ by Theorem \ref{ThB}.
The assumption $(*)$ of Theorem \ref{ThD} follows
from the successful good minimal model program for projective threefolds.
So Theorem \ref{ThD} (3) occurs, and hence $\alb_X(C_i)$ is a point
since it is $h$-periodic, noting that $C_i$ is $g$-periodic and
$g$ and $h$ are compatible. Thus, these $C_i$ are contained in
the exceptional divisor $E$ of $\alb_X$, and we may assume that
the Zariski closure $\overline{\cup_j \, C_{m_j}}$ equals $E_1$ for some irreducible component
$E_1$ of $E$ and some infinite subsequence $\{C_{m_j}\} \subset \{C_i\}$.
Thus $E_1 = D$, for $C_{m_j} \subset D$.
Since $g$ and $h$ are compatible, we have $g(E) = E$ and hence $g^s(E_1) = E_1$
for some $s > 0$. So $D = E_1$ is $g$-periodic.
This completes the proof of Theorem \ref{ThA}.

\end{document}